\pdfoutput=1
\RequirePackage{ifpdf}
\ifpdf 
\documentclass[pdftex]{sigma}
\else
\documentclass{sigma}
\fi

\newcommand{\Z}{\mathbb{Z}}
\newcommand{\Q}{\mathbb{Q}}
\newcommand{\R}{\mathbb{R}}
\newcommand{\C}{\mathbb{C}}
\newcommand{\Cstar}{{\C}^\times}
\newcommand{\Proj}{\mathbb{P}}
\newcommand{\GL}[2]{\mathrm{GL}_{#1}({#2})}

\newcommand{\orig}{\mathbf{0}}
\newcommand{\pro}[2]{{#1}({#2})}

\newcommand{\NQ}{N_\Q}
\newcommand{\V}[1]{\mathcal{V}\!\left({#1}\right)}

\newcommand{\F}[1]{\mathcal{F}\!\left({#1}\right)}
\newcommand{\bd}{\partial}
\newcommand{\Hom}{\mathrm{Hom}}
\newcommand{\dual}[1]{{{#1}^*}}
\newcommand{\intr}[1]{{{#1}^\circ}}
\newcommand{\prim}[1]{\overline{#1}}
\newcommand{\Vol}[1]{\mathrm{Vol}\!\left({#1}\right)}
\newcommand{\conv}[1]{\mathrm{conv}\!\left({#1}\right)}
\newcommand{\sconv}[1]{\mathrm{conv}\!\left\{{#1}\right\}}
\newcommand{\scone}[1]{\mathrm{cone}\!\left\{{#1}\right\}}
\newcommand{\abs}[1]{\left\vert{#1}\right\vert}
\newcommand{\width}[2]{\mathrm{width}_{#1}\!\left({#2}\right)}

\renewcommand{\min}[1]{\mathrm{min}\!\left\{{#1}\right\}}
\renewcommand{\max}[1]{\mathrm{max}\!\left\{{#1}\right\}}
\newcommand{\hmin}{{h_\mathrm{min}}}
\newcommand{\hmax}{{h_\mathrm{max}}}
\newcommand{\umin}{{u_\mathrm{min}}}
\newcommand{\umax}{{u_\mathrm{max}}}
\newcommand{\mut}{\mathrm{mut}}
\newcommand{\Newt}[1]{\mathrm{Newt}\!\left({#1}\right)}
\newcommand{\Ehr}[1]{\mathrm{Ehr}_{#1}\!\left(t\right)}
\newcommand{\Hilb}[1]{\mathrm{Hilb}_{#1}(-K_{#1})}
\newcommand{\tti}{\mathtt{i}}
\newcommand{\sfx}{\mathsf{x}}
\newcommand{\magma}{{\sc Magma}}

\begin{document}

\allowdisplaybreaks

\renewcommand{\thefootnote}{$\star$}

\renewcommand{\PaperNumber}{094}

\FirstPageHeading

\ShortArticleName{Minkowski Polynomials and Mutations}

\ArticleName{Minkowski Polynomials and Mutations\footnote{This
paper is a contribution to the Special Issue ``Mirror Symmetry and Related Topics''. The full collection is available at \href{http://www.emis.de/journals/SIGMA/mirror_symmetry.html}{http://www.emis.de/journals/SIGMA/mirror\_{}symmetry.html}}}

\Author{Mohammad AKHTAR~$^\dag$, Tom COATES~$^\dag$, Sergey GALKIN~$^\ddag$ and
  Alexander M.~KASPRZYK~$^\dag$}

\AuthorNameForHeading{M.~Akhtar, T.~Coates, S.~Galkin and A.M.~Kasprzyk}

\Address{$^\dag$~Department of Mathematics, Imperial College London,\\
\hphantom{$^\dag$}~180 Queen's Gate, London SW7~2AZ, UK}
\EmailD{\href{mailto:mohammad.akhtar03@imperial.ac.uk}{mohammad.akhtar03@imperial.ac.uk},
\href{mailto:t.coates@imperial.ac.uk}{t.coates@imperial.ac.uk},\\
\hspace*{16.8mm}\href{mailto:a.m.kasprzyk@imperial.ac.uk}{a.m.kasprzyk@imperial.ac.uk}}

\Address{$^\ddag$~Universit\"at Wien, Fakult\"at f\"ur Mathematik, Garnisongasse 3/14, 
 A-1090 Wien, Austria}
\EmailD{\href{mailto:Sergey.Galkin@phystech.edu}{Sergey.Galkin@phystech.edu}}

\ArticleDates{Received June 14, 2012, in f\/inal form December 01, 2012; Published online December 08, 2012}

\Abstract{Given \looseness=1 a Laurent polynomial $f$, one can form the period of
  $f$: this is a function of one complex variable that plays an
  important role in mirror symmetry for Fano manifolds.  Mutations are
  a particular class of birational transformations acting on Laurent
  polynomials in two variables; they preserve the period and are
  closely connected with cluster algebras. We propose a
  higher-dimensional analog of mutation acting on Laurent polynomials
  $f$ in $n$ variables. In particular we give a combinatorial
  description of mutation acting on the Newton polytope $P$ of $f$,
  and use this to establish many basic facts about mutations.
  Mutations can be understood combinatorially in terms of Minkowski
  rearrangements of slices of $P$, or in terms of piecewise-linear
  transformations acting on the dual polytope $\dual{P}$ (much like
  cluster transformations).  Mutations map Fano polytopes to Fano
  polytopes, preserve the Ehrhart series of the dual polytope, and
  preserve the period of $f$.  Finally we use our results to show
  that Minkowski polynomials, which are a~family of Laurent
  polynomials that give mirror partners to many three-dimensional Fano
  manifolds, are connected by a sequence of mutations if and only if
  they have the same period.}


\Classification{52B20; 16S34; 14J33}

\renewcommand{\thefootnote}{\arabic{footnote}}
\setcounter{footnote}{0}

\section{Introduction}\label{sec:introduction}

Given a Laurent polynomial $f \in \C\big[x^{\pm 1}_1,\ldots,x^{\pm
  1}_n\big]$, one can form the \emph{period} of~$f$
\begin{gather}
  \label{eq:period}
  \pi_f(t)=\left(\frac1{2\pi \tti}\right)^n \int_{|x_1|=\cdots
    =|x_n|=1} \frac1{1-tf(x_1,\dots, x_n)} \frac{dx_1}{x_1}\cdots
  \frac{dx_n}{x_n}, \qquad    t \in \C, \quad |t| \ll\infty .
\end{gather}
The \looseness=1 period of $f$ gives a solution to a GKZ hypergeometric
dif\/ferential system associated to the Newton polytope of $f$ (see~\cite[\S~3]{CCGGK12}). Periods of Laurent polynomials and the
associated dif\/fe\-ren\-tial systems are interesting from the point of view
of mirror symmetry, because certain Laurent polynomials arise as
mirror partners to $n$-dimensional Fano manifolds~\mbox{\cite{Ba04,BCKvS98,BCKvS00,FS, EHX,HV}}. In its most basic
form (which will suf\/f\/ice for what follows) the statement that a
Laurent polynomial~$f$ is a mirror partner for a Fano manifold~$X$
means that the Taylor expansion of the period of~$f$
\[
\pi_f(t) = \sum_{k \geq 0} c_k t^k
\]
coincides \looseness=-1 with a certain generating function for Gromov--Witten
invariants of $X$ called the \emph{quantum period} of
$X$~\cite[\S~4]{CCGGK12}. The Taylor coef\/f\/icient~$c_k$ here is the
coef\/f\/icient of the unit monomial in the Laurent polynomial~$f^k$. We
refer to the sequence $(c_k)_{k \geq 0}$ as the \emph{period sequence}
for~$f$.

We expect that if the Laurent polynomial $f$ is a mirror partner to a
Fano manifold $X$, then there is a geometric relationship between $f$
and $X$ as follows (cf.~\cite{Pr07}). Let $N_f$ be the lattice generated by the exponents
of monomials of $f$. Consider the Newton polytope $\Newt{f}$ of $f$,
and assume that the origin lies in its strict interior. Let $X_f$
denote the toric variety def\/ined by the spanning fan of $\Newt{f}$ in
$N_f\otimes\Q$; in general $X_f$ will be singular. We expect that $X_f$ admits
a smoothing with general f\/iber $X$. Note that our assumption that
$\Newt{f}$ contains the origin is not restrictive: if the origin is
outside $\Newt{f}$ then the period of $f$ is constant, and hence
cannot be the quantum period of a Fano manifold; if the origin is
contained in a~proper face of~$\Newt{f}$ then we can reduce to a
lower-dimensional situation. Note also that the lattice~$N_f$ may be a
proper sublattice of~$\Z^n$, see Example~\ref{exa:sublattice} and~\cite{BK06}. The
picture described here implies that one might expect many Laurent
polynomial mirrors for a given Fano manifold, as a smooth Fano
manifold can degenerate to many dif\/ferent singular toric varieties.

The motivating case for this paper is that of three-dimensional
\emph{Minkowski polynomials}.  These are a family of Laurent
polynomials in three variables, def\/ined in~Section~\ref{sec:Minkowski}
below, which provide mirror partners to many of the three-dimensional
Fano manifolds. The correspondence between three-dimensional Minkowski
polynomials and Fano manifolds is not one-to-one, in part because many
Minkowski polynomials give rise to the same period sequence. There are
several thousand Minkowski polynomials $f \colon (\Cstar)^3 \to \C$,
up to change of co-ordinates on~$(\Cstar)^3$, but between them these
give only 165 distinct period sequences\footnote{$98$ of these periods
  are the quantum periods for the three-dimensional Fano manifolds
  with very ample anticanonical bundle. The remaining $67$ periods are
  not the quantum period for any three-dimensional Fano manifold,
  although they may correspond to Fano
  orbifolds. See~\cite[\S~7]{CCGGK12}.}. In what follows we give a
conceptual explanation for this phenomenon. We def\/ine birational
transformations, called \emph{mutations}, that preserve periods and
show that any two Minkowski polynomials with the same period are
related by a~sequence of mutations. Our birational transformations are
higher-dimensional generalisations of the mutations def\/ined by
Galkin--Usnich~\cite{GU10}. Combining our results
in~Section~\ref{sec:Minkowski_mutations} with a~theorem of Ilten~\cite{I12}
shows that whenever~$f$ and~$g$ are Minkowski polynomials with the
same period sequence, the toric varieties~$X_f$ and~$X_g$ occur as
f\/ibers of a f\/lat family over a curve. This is consistent with our
conjectural picture, which implies that whenever Laurent polynomials
$f$ and $g$ are mirror partners for the same Fano manifold~$X$, the
toric varieties~$X_f$ and~$X_g$ are deformation equivalent.

The paper is organised as follows.  We def\/ine mutations algebraically
in~Section~\ref{sec:mutations_algebraic} and combinatorially
in~Section~\ref{sec:mutations_combinatorial}. Algebraic mutations operate on
Laurent polynomials, whereas combinatorial mutations operate on
polytopes. An algebraic mutation of a Laurent polynomial $f$ induces a
combinatorial mutation of its Newton polytope $\Newt{f}$; the converse
statement is discussed in Remark~\ref{rem:combinatorial_and_algebraic}.
We establish various basic properties of
combinatorial mutations: they send Fano polytopes to Fano polytopes
(Proposition~\ref{prop:Fano_to_Fano}); there are, up to isomorphism,
only f\/initely many mutations of a given polytope
(Proposition~\ref{prop:finitely_many}); and mutation-equivalent
polytopes have the same Hilbert series
(Proposition~\ref{prop:hilbert_preserved}).  We def\/ine Minkowski
polynomials in~Section~\ref{sec:Minkowski}, and
in~Section~\ref{sec:Minkowski_mutations} show by means of a computer search
that all Minkowski polynomials with the same period sequence are
connected by a sequence of mutations. Period sequences for Minkowski
polynomials are listed in Appendix~\href{http://www.emis.de/journals/SIGMA/2012/094/sigma12-094Appendices.pdf}{A}, and mutations connecting the
Minkowski polynomials with the same period sequence are listed in
Appendix~\href{http://www.emis.de/journals/SIGMA/2012/094/sigma12-094Appendices.pdf}{B}.

\section{Mutations}\label{sec:mutations_algebraic}

In this section we def\/ine mutations. These are a class of
birational transformations $\varphi \colon(\Cstar)^n \dashrightarrow  (\Cstar)^n$
with the property that if two Laurent polynomials $f$ and $g$ are
related by a mutation $\varphi$, so that $g = \varphi^\star f$, then the
periods of $f$ and $g$ coincide. We begin with two examples.

\begin{example}
  Consider a Laurent polynomial
  \[
  f = A(x,y) z^{-1} + B(x,y) + C(x,y) z,
  \]
  where $A$, $B$, $C$ are Laurent polynomials in $x$ and $y$. The
  pullback of $f$ along the birational transformation $(\Cstar)^3 \dashrightarrow
  (\Cstar)^3$ given by
  \begin{equation}\label{eq:mutation_example_1}
    (x,y,z) \mapsto \big(x,y, A(x,y) z\big)
  \end{equation}
  is
  \[
  g = z^{-1} + B(x,y) + A(x,y) C(x,y) z.
  \]
  We say that the Laurent polynomials $f$ and $g$ are related by the
  mutation~\eqref{eq:mutation_example_1}.
\end{example}

\begin{example}
  \label{ex:mutation_example_2}
  Consider a Laurent polynomial
  \[
  f = \sum_{i = k}^{l} C_i(x,y) z^i,
  \]
  with $k<0$ and $l>0$ where each $C_i$, $i \in \{k, k+1,\ldots,l\}$, is a Laurent
  polynomial in $x$ and $y$. Let $A$ be a Laurent polynomial in $x$
  and $y$ such that $C_i$ is divisible by $A^{-i}$ for $i \in
  \{k, k+1,\ldots,{-1}\}$. The pullback of $f$ along the birational
  transformation $(\Cstar)^3 \dashrightarrow (\Cstar)^3$ given by
  \begin{equation}\label{eq:mutation_example_2}
    (x,y,z) \mapsto \big(x,y, A(x,y) z\big)
  \end{equation}
  is
  \[
  g = \sum_{i = k}^{l} A(x,y)^i C_i(x,y) z^i.
  \]
  We say that the Laurent polynomials $f$ and $g$ are related by the
  mutation~\eqref{eq:mutation_example_2}.
\end{example}

\begin{remark}
  Note that the pullback of a Laurent polynomial along a birational
  transformation of the form~\eqref{eq:mutation_example_1}
  or~\eqref{eq:mutation_example_2} will not, in general, be a Laurent
  polynomial: the condition $A^{-i} \big| C_i$, $i \in \{k,
  k+1,\ldots,{-1}\}$, is essential.
\end{remark}

\begin{definition}\label{def:SL3_equivalence} A
  \emph{$\GL{3}{\Z}$-equivalence} is an isomorphism $(\Cstar)^3 \to
  (\Cstar)^3$ of the form
  \begin{gather*}
    (x,y,z) \mapsto \big(x^a y^b z^c, x^d y^e z^f, x^g y^h z^i\big),
\qquad
    \text{where $M := \begin{pmatrix} a & b & c \\ d & e & f \\ g & h &
        i \end{pmatrix} \in\GL{3}{\Z}$.}
  \end{gather*}
  For brevity, we write this isomorphism as $\sfx \mapsto \sfx^M$.
\end{definition}

\begin{definition}\label{def:mutation}
  A \emph{mutation} is a birational transformation $(\Cstar)^3 \dashrightarrow
  (\Cstar)^3$ given by a composition of:
  \begin{enumerate}\itemsep=0pt
  \item[1)] a $\GL{3}{\Z}$-equivalence;
  \item[2)] a birational transformation of the form~\eqref{eq:mutation_example_2}; and
  \item[3)] another $\GL{3}{\Z}$-equivalence.
  \end{enumerate}
  If $f$, $g$ are Laurent polynomials and $\varphi$ is a mutation such that
  $\varphi^\star f = g$ then we say that~$f$ and~$g$ are \emph{related by the
  mutation $\varphi$}.
\end{definition}

\begin{remark}
  One can also def\/ine mutations of Laurent polynomials in $n$
  variables, using the obvious generalisations of
  Example~\ref{ex:mutation_example_2},
  Def\/inition~\ref{def:SL3_equivalence}, and
  Def\/inition~\ref{def:mutation}.
\end{remark}

\begin{example}
  Consider the Laurent polynomial
  \[
  f = xyz+x+y+z+\frac{1}{x}+\frac{1}{xyz}.
  \]
  The Newton polytope $P$ of $f$ has two pairs of parallel facets, and
  we place one pair of them at heights $1$ and $-1$ by applying the
  $\GL{3}{\Z}$-equivalence $\sfx \mapsto \sfx^M$ with
  \[
  M =
  \begin{pmatrix}
    1 & -1 & 1 \\
    0 & 1 & -1 \\
    0 & 0 & 1
  \end{pmatrix}.
  \]
  This transforms $f$ into the Laurent polynomial
  \[
  \frac{1}{z}\left(y+\frac{y}{x}+\frac{1}{x}\right) + z\left(1+x+\frac{x}{y}\right).
  \]
  We now apply the birational transformation~\eqref{eq:mutation_example_1} with $A(x,y) =
  y+\frac{y}{x}+\frac{1}{x}$, followed by the $\GL{3}{\Z}$-equivalence
  $\sfx \mapsto \sfx^{M^{-1}}$, obtaining
  \[
  g =
  xy^2z^2+xyz+2yz^2+2z+\frac{1}{z}+\frac{1}{y}+\frac{z^2}{x}+\frac{z}{xy}.
  \]
  This shows that the Laurent polynomials $f$ and $g$ are related by
  the mutation $\varphi$, where
  \begin{gather*}
    \varphi(x,y,z)   = \left(\frac{xy^2z+yz+1}{y},
      \frac{xy^2}{xy^2z+yz+1}, \frac{z\big(xy^2z+yz+1\big)}{xy} \right).
  \end{gather*}
\end{example}

\begin{lemma}
  If the Laurent polynomials $f$ and $g$ are related by a mutation
  $\varphi$, then the periods of $f$ and $g$ coincide.
\end{lemma}

\begin{proof}
Let $\varphi \colon (\Cstar)^3 \dashrightarrow (\Cstar)^3$ be the birational transformation
\[
(x,y,z) \mapsto (x,y,A(x,y)z)
\]
from~\eqref{eq:mutation_example_2}. Since $\GL{3}{\Z}$-equivalence
preserves periods, it suf\/f\/ices to show that if $g=\varphi^\star f$
then the periods of $f$ and $g$ coincide.  Let
  \[
  Z = \big\{ (x,y,z) \in (\Cstar)^3 : A(x,y) = 0 \big \}
  \]
  and let $U = (\Cstar)^3 \setminus Z$.  Note that the restriction
  $\varphi|_U \colon U \to (\Cstar)^3$ is a morphism, and that
  \begin{equation}
    \label{eq:volume_form_is_preserved}
    (\varphi|_U)^\star \left( \frac{dx}{x} \frac{dy}{y} \frac{dz}{z}
    \right)
    =
    \frac{dx}{x} \frac{dy}{y} \frac{dz}{z}.
  \end{equation}
  Let
  \[
  C_{a,b,c} = \big\{ (x,y,z) \in (\Cstar)^3 : \text{$|x| = a$,
    $|y|=b$, $|z|=c$}\big \},
  \]
  so that the period of $f$ is
  \[
  \left(\frac1{2\pi \tti}\right)^n \int_{C_{1,1,1}}
  \frac1{1-tf(x,y,z)}
  \frac{dx}{x}
  \frac{dy}{y}
  \frac{dz}{z}.
  \]
  The amoeba of $Z$ is a proper subset of $\R^3$~\cite[Chapter~6, Corollary~1.8]{GKZ}, so there exists
  $(a,b,c)$ such that $C_{a,b,c} \subset U$.  The cycles $C_{a,b,c}$
  and $C_{a',b',c'}$ are homologous in $(\Cstar)^3$ for any non-zero
  $a$, $b$, $c$, $a'$, $b'$, $c'$.  Thus
  \begin{gather*}
    \pi_g(t)  =
    \left(\frac1{2\pi \tti}\right)^n \int_{C_{1,1,1}}
  \frac1{1-tg(x,y,z)}
  \frac{dx}{x}
  \frac{dy}{y}
  \frac{dz}{z}
    =
    \left(\frac1{2\pi \tti}\right)^n \int_{C_{a,b,c}}
  \frac1{1-tg(x,y,z)}
  \frac{dx}{x}
  \frac{dy}{y}
  \frac{dz}{z} \\
\hphantom{\pi_g(t)}{}
    =
  \left(\frac1{2\pi \tti}\right)^n \int_{\varphi(C_{a,b,c})}
  \frac1{1-tf(x,y,z)}
  \frac{dx}{x}
  \frac{dy}{y}
  \frac{dz}{z} \quad \text{by the change of variable formula.}
  \end{gather*}
  Now the homology class $\big[\varphi(C_{a,b,c})\big] \in H_3\big(
  (\Cstar)^3,\Z\big)$ is equal to $k \big[C_{1,1,1} \big]$ for some
  integer $k$, since $H_3\big( (\Cstar)^3,\Z\big)$ is freely generated
  by $\big[C_{1,1,1} \big]$, and from~  \eqref{eq:volume_form_is_preserved} and the change of variable
  formula we see that $k=1$.  It follows that
\begin{gather*}
  \pi_g(t) =
  \left(\frac1{2\pi \tti}\right)^n \int_{C_{1,1,1}}
  \frac1{1-tf(x,y,z)}
  \frac{dx}{x}
  \frac{dy}{y}
  \frac{dz}{z}  = \pi_f (t).\tag*{\qed}
\end{gather*}
\renewcommand{\qed}{}
\end{proof}

\begin{example}
Consider the two Laurent polynomials
\begin{gather*}
f_1 =x+\frac{2x}{y}+\frac{x}{y^2}+y+z+\frac{1}{z}+\frac{z}{y}+\frac{4}{y}+\frac{1}{yz}+\frac{yz}{x}+\frac{2y}{x}\\
\hphantom{f_1 =}{}
+\frac{y}{xz}
+\frac{2z}{x}+\frac{5}{x}+\frac{2}{xz}+\frac{yz}{x^2}+\frac{2y}{x^2}+\frac{y}{x^2z},\\
f_2 =x+\frac{2x}{y}+\frac{x}{y^2}+y+z+\frac{1}{z}+\frac{z}{y}+\frac{3}{y}+\frac{1}{yz}+\frac{yz}{x}+\frac{3y}{x}\\
\hphantom{f_2 =}{}
+\frac{y}{xz}
+\frac{2z}{x}+\frac{4}{x}+\frac{2}{xz}+\frac{yz}{x^2}+\frac{2y}{x^2}+\frac{y}{x^2z}.
\end{gather*}
Since $f_1$ and $f_2$ have the same period sequence
\[
\pi_{f_1}(t)=\pi_{f_2}(t)=1+28t^2+216t^3+3516t^4+49680t^5+\cdots,
\]
and since $\Newt{f_1}=\Newt{f_2}$, it is tempting to assume that there
is some $\GL{3}{\Z}$-equivalence that preserves the Newton polytope
and sends $f_1$ to $f_2$. This is not the case. However, there does
exist a birational map sending $f_1$ to $f_2$. This is a composition
of mutations
\[
f_1\stackrel{\varphi}{\longrightarrow}f\stackrel{\psi}{\longrightarrow}f_2
\]
factoring through
\[
f=xz^2+2xz+x+yz+y+3z+\frac{2}{z}+\frac{z}{y}+\frac{1}{y}+\frac{y}{x}+\frac{y}{xz}+\frac{2}{x}+\frac{2}{xz}+\frac{1}{xz^2}+\frac{1}{xy}+\frac{1}{xyz}.
\]
The maps $\varphi$, $\psi$ and their inverses are given by
\begin{gather*}
\varphi: \ (x,y,z) \mapsto\left(\frac{z(xyz+(y+1)^2)}{y}, \frac{xyz+(y+1)^2}{xy}, y\right),\\
\psi: \ (x,y,z) \mapsto\left(\frac{(x+yz+y)(xz+yz+y)}{y^2z(x+y)}, \frac{1}{z}, \frac{y}{x}\right),\\
\varphi^{-1}: \ (x,y,z) \mapsto\left(\frac{xz+y(z+1)^2}{y^2z}, z, \frac{xyz}{xz+y(z+1)^2}\right),\\
\psi^{-1}: \ (x,y,z) \mapsto\left(\frac{(yz+z+1)(yz+y+z)}{xyz^2(z+1)}, \frac{(yz+z+1)(yz+y+z)}{xyz(z+1)}, \frac{1}{y}\right).
\end{gather*}
Set $P:=\Newt{f_1}=\Newt{f_2}$ and $Q:=\Newt{f}$. The polytopes $P$
and $Q$ are ref\/lexive, but since $\Vol{P}=32$ and $\Vol{Q}=28$, $P$
and $Q$ are not isomorphic. However, as predicted by
Proposition~\ref{prop:hilbert_preserved} below, the Ehrhart series
$\Ehr{\dual{P}}$ and $\Ehr{\dual{Q}}$ are equal: in other words, the
Hilbert series $\Hilb{X_{f_i}}$ and $\Hilb{X_f}$ agree.
\end{example}

\section{Combinatorial mutations}\label{sec:mutations_combinatorial}

The mutations of a Laurent polynomial $f$ def\/ined in~Section~\ref{sec:mutations_algebraic} induce transformations of the Newton polytope of~$f$. In this section we give a combinatorial and coordinate-free def\/inition of these transformations, which we call \emph{combinatorial mutations}, in terms of the Newton polytope alone. We then establish some basic properties of combinatorial mutations. Let us begin by f\/ixing our notation. Let $N$ be an $n$-dimensional lattice and let $P\subset\NQ:=N\otimes\Q$ be a convex lattice polytope such that
\begin{enumerate}\itemsep=0pt
\item[1)] $P$ is of maximum dimension, $\dim{P}=n$;
\item[2)] 
the origin lies in the strict interior of $P$, $\orig\in\intr{P}$;
\item[3)] the vertices $\V{P}\subset N$ of $P$ are primitive lattice points.
\end{enumerate}
We call such a polytope \emph{Fano}.

Given any lattice polytope $P\subset\NQ$, the \emph{dual} polyhedron $\dual{P}\subset M_\Q$, where $M:=\Hom(N,\Z)$, is def\/ined by
\[
\dual{P}:=\{u\in M_\Q\mid\pro{u}{v}\geq -1\text{ for all }v\in P\}.
\]
Condition~2 ensures that, when $P$ is Fano, $\dual{P}$
is a polytope. When $\dual{P}$ is a \emph{lattice} polytope, we say
that $P$ is \emph{reflexive}. Low-dimensional ref\/lexive polytopes have
been classif\/ied~\cite{KS98b,KS00}: up to the action of $\GL{n}{\Z}$
there are $16$ ref\/lexive polytopes in dimension two; $4,319$ in
dimension three; and $473,800,776$ in dimension four. A Fano polytope
$P\subset\NQ$ is called \emph{canonical} if $\intr{P}\cap
N=\{\orig\}$. In two dimensions, the ref\/lexive polytopes and canonical
polytopes coincide. In general every ref\/lexive polytope is canonical,
although the converse is not true: there are $674,688$ canonical
polytopes in dimension three~\cite{Kas08a}.

\begin{definition}\label{def:width_vector}
  Let $w\in M$ be a primitive lattice vector, and let $P\subset\NQ$ be
  a lattice polytope. Set
  \begin{gather*}
    \hmin:=\min{\pro{w}{v}\mid v\in P}, \qquad
    \hmax:=\max{\pro{w}{v}\mid v\in P}
  \end{gather*}
  We def\/ine the \emph{width} of $P$ with respect to $w$ to be the
  positive integer
\[
\width{w}{P}:=\hmax - \hmin.
\]
  If $\width{w}{P}=l$ then we refer to $w$ as a \emph{width $l$
    vector} for $P$. We say that a lattice point $v\in N$ (resp.~a
  subset $F\subset\NQ$) is at \emph{height $m$} with respect to $w$ if
  $\pro{w}{v}=m$ (resp.~if $\pro{w}{F}=\{m\}$).
\end{definition}

If $\orig\in\intr{P}$ then $\hmin<0$ and $\hmax>0$, hence $w$ must
have width at least two. If $P$ is a~ref\/lexive polytope then for any
$w\in\V{\dual{P}}$ there exists a facet $F\in\F{P}$ at height $-1$;
this is a well-known characterisation of ref\/lexive
polytopes~\cite{Bat94}.  For each height $h\in\Z$, $w$ def\/ines a~hyperplane $H_{w,h}:=\{x\in\NQ\mid\pro{w}{x}=h\}$. Let
\[
w_h(P):=\conv{H_{w,h}\cap P\cap N}.
\]
By def\/inition, $w_\hmin(P)=H_{w,\hmin}\cap P$ and $w_\hmax(P)=H_{w,\hmax}\cap P$ are faces of $P$, hence $\V{w_\hmin(P)}\subseteq\V{P}$ and $\V{w_\hmax(P)}\subseteq\V{P}$. Furthermore, the face $w_\hmin(P)$ is a facet of $P$ if and only if $w=\prim{u}$ for some vertex $u\in\V{\dual{P}}$, where $\prim{u}$ denotes the unique primitive lattice point on the ray from $\orig$ through $u$. Similarly, $w_\hmax(P)$ is a facet if and only if $-w=\prim{u}$ for some $u\in\V{\dual{P}}$.

\begin{definition}\label{def:empty_sum}
Recall that the Minkowski sum of two polytopes $Q,R\subset\NQ$, is
\[
Q+R:=\{q+r\mid q\in Q,r\in R\}.
\]
Henceforth we adopt the convention that $Q+\varnothing:=\varnothing$ for any polytope $Q$.
\end{definition}

\begin{definition}\label{def:comb_mut}
Suppose that there exists a lattice polytope $F\subset\NQ$ with $\pro{w}{F}=0$, such that for every height $\hmin\leq h<0$ there exists a possibly-empty lattice polytope $G_h\subset\NQ$ satisfying $H_{w,h}\cap\V{P}\subseteq G_h+(-h)F\subseteq w_h(P)$. We call such an $F$ a \emph{factor} for $P$ with respect to $w$. We def\/ine the \emph{combinatorial mutation} given by width vector~$w$, factor~$F$, and polytopes $\{G_h\}$ to be the convex lattice polytope
\[
\mut_w(P,F;\{G_h\}):=\conv{\bigcup_{h=\hmin}^{-1}G_h\cup\bigcup_{h=0}^{\hmax}(w_h(P)+hF)}\subset\NQ.
\]
\end{definition}

Notice that one need only consider factors up to translation, since for any $v\in N$ such that $\pro{w}{v}=0$ we have $\mut_w(P,v+F;\{G_h+hv\})\cong\mut_w(P,F;\{G_h\})$. In particular, if $w_\hmin(P)$ is zero-dimensional then combinatorial mutations with width vector $w$ leave $P$ unchanged.

\begin{example}
  \label{ex:algebraic_gives_combinatorial}
  Consider the situation of Example~\ref{ex:mutation_example_2}, so
  that $f$ is the Laurent polynomial
  \[
  f = \sum_{i = k}^{l} C_i(x,y) z^i,
  \]
  with $k<0$ and $l>0$, $A$ is a Laurent polynomial in $x$ and $y$
  such that
\begin{gather*}
 A^{-i} \big| C_i \qquad \text{for $i \in \{k,k+1,\ldots,{-1}\}$},
\end{gather*}
      $\varphi$ is the birational transformation $(x,y,z) \mapsto
      \big(x,y, A(x,y) z\big)$, and
      \[
      g = \sum_{i = k}^{l} A(x,y)^i C_i(x,y) z^i.
      \]
      The Laurent polynomials $f$ and $g$ are related by the algebraic
      mutation $\varphi$.  This algebraic mutation induces a combinatorial
      mutation of $P = \Newt{f}$ with $w = (0,0,1)$, $F = \Newt{A}$,
      $\hmin = k$, $\hmax = l$,
\begin{gather*}
  G_h = \Newt{\frac{C_{h}}{A^{-h}}}, \qquad \text{where\footnotemark\ $\hmin \leq h \leq -1$},
  \end{gather*}
  \footnotetext{The Newton polytope of the zero
    polynomial is $\varnothing$.}
  and $\mut_w(P,F;\{G_h\}) = \Newt{g}$.
\end{example}

\begin{remark}
  \label{rem:combinatorial_and_algebraic}
  Given a Laurent polynomial $f_1$ with Newton polytope $P$, there may
  exist combinatorial mutations of $P$ that do not arise from any
  algebraic mutation of $f_1$.  Given a combinatorial mutation of $P$,
  however, there exists a Laurent polynomial $f_2$ with $\Newt{f_2} =
  P$ such that the combinatorial mutation arises from an algebraic
  mutation of $f_2$.  See Example~\ref{ex:combinatorial_and_algebraic}
  below.
\end{remark}

\begin{lemma}\label{lem:mutation_inverse}
Let $Q:=\mut_w(P,F;\{G_h\})$ be a combinatorial mutation of $P\subset\NQ$. Then $\mut_{-w}(Q,F;\{w_h(P)\})$ is a combinatorial mutation of $Q$ equal to $P$.
\end{lemma}
\begin{proof}
Let $P':=\mut_{-w}(Q,F;\{w_h(P)\})$. Clearly this is well def\/ined. Let $v\in\V{P}$ be a~vertex of $P$ with height $h$. If $h\ge 0$ then $P'\supseteq w_h(P)\supseteq H_{w,h}\cap\V{P}$, so $v\in P'$. If $h<0$ then $P'\supseteq G_h+(-h)F\supseteq H_{w,h}\cap\V{P}$. Hence $P\subseteq P'$. Conversely let $v\in\V{P'}$ and set $h:=w(v)$. If $h\geq 0$ then $v\in w_h(P)$, so $v\in P$. If $h<0$ then $v\in G_h+(-h)F\subseteq w_h(P)$, so again $v\in P$. Hence $P=P'$.
\end{proof}

\begin{lemma}\label{lem:vertices_mutations}
Let $Q:=\mut_w(P,F;\{G_h\})$ be a combinatorial mutation of $P\subset\NQ$. Then
\begin{gather*}
\V{Q} \subseteq\{v_P+\pro{w}{v_P}v_F\mid v_P\in w_h(P)\cap N,v_F\in\V{F}\}, \qquad \text{and}\\
\V{P} \subseteq\{v_Q-w(v_Q)v_F\mid v_Q\in w_h(Q)\cap N,v_F\in\V{F}\}.
\end{gather*}
\end{lemma}

\begin{proof}
Pick any vertex $v\in\V{Q}$, and set $h:=w(v)$. First consider the case when $h\geq 0$. Notice that $v\in\V{w_h(P)+hF}$, since otherwise $v$ lies in the convex hull of two (not necessarily lattice) points $v_1,v_2\in Q$, and as such could not be a vertex. In general the vertices of a Minkowski sum are contained in the sum of the vertices of the summands, hence there exist lattice points $v_P\in\V{w_h(P)}$ and $v_F\in\V{F}$ such that $v=v_P+hv_F$. If $h<0$ then $v\in\V{G_h}$. In particular, $v+(-h)v_F\subseteq w_h(P)$ for all $v_F\in F$, so there exist lattice points $v_P\in w_h(P)$ and $v_F\in\V{F}$ such that $v=v_P+hv_F$. Hence we have the f\/irst equation in the statement.

The second equation follows from the f\/irst by considering the inverse combinatorial mutation $\mut_{-w}(Q,F;\{w_h(P)\})$.
\end{proof}

\begin{proposition}\label{prop:mutations_are_unique}
Let $P\subset\NQ$ be a convex lattice polytope with width vector $w$ and factor $F$. Then $\mut_w(P,F;\{G_h\})=\mut_w(P,F;\{G'_h\})$ for any two combinatorial mutations of $P$.
\end{proposition}
\begin{proof}
Set $Q:=\mut_w(P,F;\{G_h\})$ and $Q':=\mut_w(P,F;\{G'_h\})$, and suppose that $Q\neq Q'$. Then (possibly after exchanging $Q$ and $Q'$) there exists some vertex $q\in\V{Q'}$ such that $q\not\in Q$, $\pro{w}{q}<0$. In particular, there exists a supporting hyperplane $H_{u,l}$ of $Q$ separating $Q$ from $q$; i.e.~$\pro{u}{x}\leq l$ for all $x\in Q$, and $\pro{u}{q}>l$.

For any $v\in\V{P}$ there exists $v_Q\in w_{\pro{w}{v}}(Q)\cap N$, $v_F\in\V{F}$ such that $v=v_P+\pro{w}{v_Q}v_F$. Hence
\[
\pro{u}{v}=\pro{u}{v_Q}-\pro{w}{v_Q}\pro{u}{v_F}\leq
\begin{cases}
l-\pro{w}{v_Q}\umin,&\text{if }\pro{w}{v_Q}\geq 0,\\
l-\pro{w}{v_Q}\umax,&\text{if }\pro{w}{v_Q}<0,
\end{cases}
\]
where $\umin:=\min{\pro{u}{v_F}\mid v_F\in\V{F}}$, and $\umax:=\max{\pro{u}{v_F}\mid v_F\in\V{F}}$.

Since $\mut_{-w}(Q',F;\{w_h(P)\})=P$, so $q-\pro{w}{q}F\subseteq P$. By def\/inition there exists some $v_F\in\V{F}$ such that $\pro{u}{v_F}=\umax$, hence
\[
\pro{u}{q-\pro{w}{q}v_F}=\pro{u}{q}-\pro{w}{q}\umax>l-\pro{w}{q}\umax,
\]
a contradiction. Hence $Q=Q'$.
\end{proof}

In light of Proposition~\ref{prop:mutations_are_unique}, we simply write $\mut_w(P,F)$ for a mutation $\mut_w(P,F;\{G_h\})$ of $P$.

\begin{example} \label{ex:combinatorial_and_algebraic}\sloppy
Consider the Laurent polynomial
\[
f=\frac{z^2}{y}+2z^2+yz^2+\frac{2xz^2}{y}+2xz^2+\frac{x^2z^2}{y}+\frac{1}{x^2}+\frac{x^4}{y^2}+y^2+\frac{1}{z}.
\]
The corresponding Newton polytope $P$ has vertices $\{(0,1,2),
(0,-1,2), (2,-1,2), (0,2,0)$, $(4,-2,0), (-2,0,0),
(0,0,-1)\}$. Note that the sublattice generated by the non-zero
coef\/f\/icients has index one.

Set $w=(0,0,-1)\in M$. The height $-2$ slice $w_{-2}(P)$
Minkowski-factorizes into two empty triangles, whereas the height $-1$
slice $w_{-1}(P)=\sconv{(0,1,1),(-1,0,1),(1,-1,1),(2,-1,1)}$ is
indecomposable. Since there are no vertices of $P$ at height $-1$,
Def\/inition~\ref{def:comb_mut} allows us to take
$G_{-1}=\varnothing$. This gives a combinatorial mutation to
\[
Q=\sconv{(0,-1,2), (4,-2,0), (0,2,0), (-2,0,0), (1,0,-1), (0,1,-1), (0,0,-1)}.
\]
This corresponds to the algebraic mutation
$\varphi:(x,y,z)\mapsto(x,y,z/(x+y+1))$ sending $f$ to
\[
\varphi^*f=\frac{z^2}{y}+\frac{x^4}{y^2}+y^2+\frac{1}{x^2}+\frac{x}{z}+\frac{y}{z}+\frac{1}{z}.
\]
Set $g = f + yz$.  Then $\Newt{g} = P$ but, since $g$ has a non-zero
coef\/f\/icient on the slice $w_{-1}(P)$, the combinatorial
mutation described above does not arise from any algebraic mutation of~$g$.
\end{example}

\begin{example}
The weighted projective spaces $\Proj(1,1,1,3)$ and $\Proj(1,1,4,6)$ are known to have the largest degree $-K^3=72$ amongst all canonical \emph{toric} Fano threefolds~\cite[Theorem~3.6]{Kas08a} and amongst all \emph{Gorenstein} canonical Fano threefolds~\cite{Pro05}. We shall show that they are connected by a width three combinatorial mutation.

Let $P:=\sconv{(1,0,0),(0,1,0),(0,0,1),(-1,-1,-3)}\subset\NQ$ be the simplex associated with $\Proj(1,1,1,3)$. The primitive vector $(-1,2,0)\in M$ is a width three vector on $P$, with $w_{-1}(P)$ equal to the edge $\sconv{(-1,-1,-3),(1,0,0)}$, and $w_{2}(P)$ given by the vertex $(0,1,0)$. Let $F:=\sconv{(0,0,0),(2,1,3)}$, and consider the mutation $\mut_w(P,F)$. This has vertices
\[
\{(-1,-1,-3),(0,0,1),(0,1,0),(4,3,6)\},
\]
and is the simplex associated with $\Proj(1,1,4,6)$.
\end{example}

\begin{proposition}\label{prop:Fano_to_Fano}
Let $P\subset\NQ$ be a lattice polytope. The combinatorial mutation $\mut_w(P,F)$ is a Fano polytope if and only if $P$ is a Fano polytope.
\end{proposition}
\begin{proof}
Begin by assuming that $P$ is Fano, and set $Q:=\mut_w(P,F)$. Let $v\in\V{Q}$, and def\/ine $h:=\pro{w}{v}$ to be the height of $v$ with respect to $w$. If $h\geq 0$ then $v=v_P+hv_F$ for some $v_P\in\V{P}$, $v_F\in\V{F}$. But $v_P$ is primitive by assumption, hence $v$ is primitive.

If $h<0$ then $v\in\V{G_h}$. Without loss of generality we are free to take $G_h$ equal to the smallest polytope such that $G_h+(-h)F\supseteq H_{w,h}\cap\V{P}$. Suppose that $v$ is not primitive. Then for any $v_F\in\V{F}$, $v+(-h)v_F$ is not primitive, hence $v+(-h)v_F$ is not a point of $H_{w,h}\cap\V{P}$. But this implies that we can take $G'_h:=\conv{G_h\cap N\setminus\{v\}}\subsetneq G_h$ with $G'_h+(-h)F\supseteq H_{w,h}\cap\V{P}$, contradicting our choice of $G_h$. Hence $v$ is primitive.

Finally, the ``if and only if'' follows by considering the inverse combinatorial mutation.
\end{proof}

\begin{corollary}\label{cor:canonical_width_2}
Let $P\subset\NQ$ be a lattice polytope, and let $w$ be a width two vector. The combinatorial mutation $\mut_w(P,F)$ is a canonical polytope if and only if $P$ is a canonical polytope.
\end{corollary}
\begin{proof}
Begin by assuming that $P$ is a canonical polytope, and set $Q:=\mut_w(P,F)$. Since~$w$ is of width two, we need only show that $H_{w,0}\cap\V{P}=H_{w,0}\cap\V{Q}$; it follows that $Q$ is canonical. But $H_{w,0}\cap\V{Q}\subseteq H_{w,0}\cap\V{P}\subseteq H_{w,0}\cap\V{Q}$, where the second inclusion follows by considering the inverse combinatorial mutation. Once more, the ``if and only if'' follows by exchanging the roles of $P$ and $Q$ via the inverse combinatorial mutation.
\end{proof}

\begin{proposition}\label{prop:finitely_many}
Let $P\subset\NQ$ be a lattice polytope, $\orig\in\intr{P}$. Up to isomorphism, there are only finitely many combinatorial mutations $\mut_w(P,F)$.
\end{proposition}

\begin{proof}
  For f\/ixed width vector $w$ there are clearly only f\/initely many
  factors $F$, and so only f\/initely many combinatorial mutations
  $\mut_w(P,F)$. Assume that $\orig\in\intr{P}$ and consider the
  spanning fan $\Delta$ in $M_\Q$ for the dual polytope
  $\dual{P}$. Since $\orig\in\intr{(\dual{P})}$, any width vector $w$
  lies in a cone $\sigma\in\Delta$. If we insist that $\dim\sigma$ is
  as small as possible, then $\sigma$ is uniquely determined. If
  $\dim\sigma=\dim P$ then $w_\hmin(P)$ is zero-dimensional, hence
  $\mut_w(P,F)\cong P$. So we may insist that
  $\sigma\in\Delta^{(n-1)}$ is of codimension at least one.

  By def\/inition of duality, $w\in\bd(-\hmin\dual{P})$. Consider the
  corresponding face $F:=w_\hmin(P)$. Let $l_F\in\Z_{>0}$ be the
  largest integer such that there exist lattice polytopes
  $A,B\subset\NQ$, $A\neq\varnothing$, with $F=l_FA+B$. Then, for a
  factor to exist, $l\leq l_F$. Def\/ine $l_P:=\max{l_F\mid\text{$F$ is
      a face of $P$}}$. Then $w$ is a primitive vector in
  $\Delta^{(n-1)}\cap l_P\dual{P}\cap M$, and the right hand side is a
  f\/inite set that depends only on $P$.
\end{proof}

\begin{remark}
When $\dim P=2$, $\orig\in\intr{P}$, we see that $\mut_w(P,F)$ is a non-trivial combinatorial mutation only if $w\in\{\prim{u}\mid u\in\V{\dual{P}}\}$. In particular, there are at most $\abs{\V{P}}$ choices for $w$, and at most $\abs{\bd P\cap N}$ distinct non-trivial combinatorial mutations.
\end{remark}

Let $Q\subset M_\Q$ be a rational polytope, and let $r\in\Z_{>0}$ be the smallest positive integer such that $rQ$ is a lattice polytope. In general there exists a quasi-polynomial, called the \emph{Ehrhart quasi-polynomial}, $L_Q:\Z\rightarrow\Z$ of degree $\dim Q$ such that $L_Q(m)=\abs{mQ\cap M}$. The corresponding generating function $\Ehr{Q}:=\sum\limits_{m\geq 0}L_Q(m)t^m$ is called the \emph{Ehrhart series} of $Q$, and can be written as a rational function~\cite{Stan80}
\[
\Ehr{Q}=\frac{\delta_0+\delta_1t+\ldots+\delta_{r(n+1)-1}t^{r(n+1)-1}}{(1-t^r)^{n+1}}
\]
with non-negative coef\/f\/icients. We call $(\delta_0,\delta_1,\ldots,\delta_{r(n+1)-1})$ the \emph{$\delta$-vector} of $Q$. In particular, the $\delta$-vector is palindromic if and only if $\dual{Q}\subset\NQ$ is a lattice polytope~\cite{FK08}.

\begin{proposition}\label{prop:hilbert_preserved}
Let $P\subset\NQ$ be a lattice polytope, $\orig\in\intr{P}$, and let $Q:=\mut_w(P,F)$ be a~combinatorial mutation of $P$. Then $\Ehr{\dual{P}}=\Ehr{\dual{Q}}$.
\end{proposition}
\begin{proof}
Since $\dual{P}$ is the dual of a lattice polytope, for any point $u\in M$ there exists a non-negative integer $k\in\Z_{\geq 0}$ such that $u\in\bd(k\dual{P})$. By def\/inition of duality, $H_{u,-k}:=\{\pro{u}{v}=-k\mid v\in\NQ\}$ is a supporting hyperplane for $P$. We begin by showing that the map $u\mapsto u-\umin w$, where $\umin:=\min{\pro{u}{v_F}\mid v_F\in\V{F}}$, gives a supporting hyperplane $H_{u-\umin w,-k}$ for $Q$, hence $u-\umin w\in\bd(k\dual{Q})$.

Let $v\in\V{Q}$. If $\pro{w}{v}\geq 0$, we can write $v=v_P+\pro{w}{v_P}v_F$ for some $v_P\in\V{P}$, $v_F\in\V{F}$. In particular,
\[
\pro{(u-\umin w)}{v}=\pro{u}{v_P}+\pro{w}{v_P}(\pro{u}{v_F}-\umin)\geq -k.
\]
Now suppose that $\pro{w}{v}<0$. For any $v_F\in\V{F}$ we have that $v-\pro{w}{v}u_F\in P$, hence $\pro{u}{v-\pro{w}{v}v_F}\geq -k$. In particular, $\pro{u}{v}\geq -k+\umin\pro{w}{v}$, so
\[
\pro{(u-\umin w)}{v}=\pro{u}{v}-\umin\pro{w}{v}\geq -k+\umin\pro{w}{v}-\umin\pro{w}{v}=-k.
\]

Since $H_{u,-k}$ is a supporting hyperplane for $P$, there exists some $v_P\in\V{P}$ such that $\pro{u}{v_P}=-k$. If $\pro{w}{v_P}\geq 0$ then for any $v_F\in\V{F}$, $v_P+\pro{w}{v_P}v_F\in Q$. Picking $v_F$ such that $\pro{u}{v_F}=\umin$, we obtain $\pro{(u-\umin w)}{v_P+\pro{w}{v_P}v_F}=-k$. If $\pro{w}{v_P}<0$ then there exists some $v_Q\in\V{Q}$, $v_F\in\V{F}$ such that $v_P=v_Q-\pro{w}{v_P}v_F$. Suppose that $\pro{u}{v_F}>\umin$, and let $v_F'\in\V{F}$ be such that $\pro{u}{v_F'}=\umin$. Then $\pro{u}{v_Q-\pro{w}{v_P}v_F'}<\pro{u}{v_Q-\pro{w}{v_P}v_F}=-k$. But $v_Q-\pro{w}{v_P}v_F'\in P$, which is a contradiction. Hence $\pro{u}{v_F}=\umin$, and so $\pro{(u-\umin w)}{v_Q}=-k$. Thus $H_{u-\umin w,-k}$ is a supporting hyperplane for $Q$.

Finally, suppose that $u,u'\in\bd(k\dual{P})$ are lattice points such that $u-\umin w=u'-u'_\mathrm{min}w$. Since for any $v_F\in\V{F}$, $\pro{(u-\umin w)}{v_F}=\pro{u}{v_F}$, we see that $\umin=u'_\mathrm{min}$ and so $u=u'$. Hence $\abs{\bd(k\dual{P})}\leq\abs{\bd(k\dual{Q})}$. By considering the inverse combinatorial mutation $\mut_{-w}(Q,F)=P$, we have equality. The result follows.
\end{proof}

\begin{example}\label{exa:sublattice}
Consider the Laurent polynomials
\begin{gather*}
f_1   =xyz^2u+x+y+z+\frac{1}{yz}+\frac{1}{x^2yz^2u},
\qquad
f_2  =xyz^2u^3+x+y+z+\frac{1}{yz}+\frac{1}{x^2yz^2u^3}.
\end{gather*}
 The Newton polytopes $P_i:=\Newt{f_i}$ are both four-dimensional
ref\/lexive polytopes in $\NQ$, and the $\delta$-vectors of the dual polytopes
are, respectively,
\begin{gather*}
\delta_1  = (1,95,294,95,1),\qquad
\delta_2  = (1,29,102,29,1).
\end{gather*}
Proposition~\ref{prop:hilbert_preserved} thus implies that there is no
sequence of mutations connecting~$f_1$ and~$f_2$. However
\[
\pi_{f_1}(t)=\pi_{f_2}(t)=1+12t^3+900t^6+94080t^9+11988900t^{12}+\cdots,
\]
which we recognise as the period sequence for $\Proj^2\times\Proj^2$.
As in~Section~\ref{sec:introduction}, let $N_{f_i}$ be the sublattice of $N$
generated by the exponents of monomials of $f_i$, and let $X_{f_i}$ be
the toric variety def\/ined by the spanning fan of $P_i$ in
$N_{f_i}\otimes\Q$. The lattice $N_{f_1}$ is equal to $N$, and $X_{f_1}$
is isomorphic to $\Proj^2 \times \Proj^2$. The lattice $N_{f_2}$ is an
index-three sublattice of $N$, but the restriction of $P_2$ to $N_{f_2}$
is isomorphic to $P_1$ and hence $X_{f_2}$ is also isomorphic to
$\Proj^2 \times \Proj^2$. Consider now the toric variety
$\widetilde{X}_2$ def\/ined by the spanning fan of $P_2$ in $\NQ$. This
is a three-fold cover of $X_2$, and is not deformation-equivalent to
$\Proj^2 \times \Proj^2$: $\widetilde{X}_2$ has Hilbert
$\delta$-vector $\delta_2$, whereas $\Proj^2 \times \Proj^2$ has
Hilbert $\delta$-vector $\delta_1$. As this example illustrates,
sublattices are invisible to period sequences.
\end{example}

Let $Q:=\mut_w(P,F)$ be a combinatorial mutation of $P$, and let $r_P$ denote the smallest dilation $k$ of $\dual{P}$ such that $k\dual{P}$ is integral. When $P$ is a Fano polytope, we call $r_P$ the \emph{Gorenstein index}. The minimum period common to the cyclic coef\/f\/icients of $L_\dual{P}$ divides $r_P$; when the period does not equal $r_P$ we have a phenomena known as \emph{quasi-period collapse}~\cite{BSW08,HM08}. In general $r_P\ne r_Q$ but, by Proposition~\ref{prop:hilbert_preserved}, $\dual{P}$ and $\dual{Q}$ are Ehrhart equivalent. Hence we have examples of quasi-period collapse. At its most extreme, when $P$ is ref\/lexive $r_P=1$ and so the period of~$L_\dual{Q}$ is one. Combinatorial mutations give a systematic way of producing families of rational polytopes that exhibit this behaviour.

\begin{corollary}\label{cor:dim3_reflexive}
A three-dimensional canonical polytope $P\subset\NQ$ is combinatorial mutation equivalent to a reflexive polytope if and only if $P$ is reflexive.
\end{corollary}
\begin{proof}
By inspecting the classif\/ication of canonical polytopes~\cite{Kas08a} we see that the period of~$L_\dual{P}$ equals~$r_P$ for each~$P$.
\end{proof}

The following example demonstrates that Corollary~\ref{cor:dim3_reflexive} does not hold in higher dimensions.

\begin{example}
For any even dimension $n=2k\ge 4$, def\/ine $P_k$ to be the polytope with $n+2$ vertices:
\begin{gather*}
P_k:=\mathrm{conv}\{ (2,2,1,1,\ldots,1),(2,1,2,1,\ldots,1),\\
\hphantom{P_k:=\mathrm{conv}\{} (0,1,\ldots,0),\ldots,(0,0,\ldots,1),(-1,-1,\ldots,-1)\}.
\end{gather*}
This is a canonical polytope, but is not ref\/lexive: the facet def\/ined by all but the f\/inal vertex is supported by the primitive vector $(n-1,-2,\ldots,-2)\in M$ at height $-2$. Let $w:=(1,-1,-1,0,\ldots,0)\in M$ be a width two vector on~$P_k$. The face $w_{-1}(P_k)$ is two-dimensional; this is the empty square, and Minkowski factors into two line segments. The corresponding combinatorial mutation gives the standard polytope for $\Proj^n$. By Proposition~\ref{prop:hilbert_preserved},
\[
\Ehr{\dual{P_k}}=\Hilb{\Proj^n}.
\]
Since $r_{P_k}=2$, we have quasi-period collapse.
\end{example}

The map $u\mapsto u-\umin w$ on $M$ described in the proof of
Proposition~\ref{prop:hilbert_preserved} is a piecewise-linear
transformation analogous to a cluster transformation. Let $\Delta_F$
be the normal fan of $F$ in $M_\Q$. Notice that the normal fan is
well-def\/ined up to translation of $F$. Furthermore, the cones of
$\Delta_F$ are not strictly convex: they each contain the subspace $\Z
w$. Let $\sigma$ be a maximum-dimensional cone of $\Delta_F$. Then
there exists $M_\sigma \in \GL{n}{\Z}$ such that, for any
$u\in-\sigma$, the map $u\mapsto u-\umin w$ is equal to the map $u
\mapsto u M_\sigma$.

We make this explicit. Suppose without loss of generality that
$w=(0,\ldots,0,1)\in M$. Then a maximum-dimensional cone
$\sigma\in\Delta_F$ corresponds to some vertex
$v_\sigma=(v_{\sigma,1},\ldots,v_{\sigma,n-1},0)\in\V{F}$, and
\[
 M_\sigma := \begin{pmatrix}
1&0&\ldots&0&-v_{\sigma,1}\\
0&1&\ldots&0&-v_{\sigma,2}\\
\vdots&\vdots&&\vdots&\vdots\\
0&0&\ldots&1&-v_{\sigma,n-1}\\
0&0&\ldots&0&1
\end{pmatrix}.
\]

\begin{example}
  \label{ex:piecewise_GL3}
  Let $w:=(0,0,1)\in M$ and let
  $F:=\sconv{(0,0,0),(1,0,0),(0,1,0)}\subset\NQ$. We have:
\[
M_\sigma= \begin{cases}
\begin{pmatrix}1&0&0\\0&1&0\\0&0&1\end{pmatrix}&
\text{if }\sigma=\scone{(-1,0,0),(0,-1,0),(0,0,\pm 1)},\\
\begin{pmatrix}1&0&-1\\0&1&0\\0&0&1\end{pmatrix}&
\text{if }\sigma=\scone{(1,1,0),(0,-1,0),(0,0,\pm 1)},\\
\begin{pmatrix}1&0&0\\0&1&-1\\0&0&1\end{pmatrix}&
\text{if }\sigma=\scone{(-1,0,0),(1,1,0),(0,0,\pm 1)}.
\end{cases}
\]
\end{example}

\section{Three-dimensional Minkowski polynomials}\label{sec:Minkowski}

Three-dimensional Minkowski polynomials are a family of Laurent
polynomials in three variables with Newton polytopes given by
three-dimensional ref\/lexive polytopes. We consider Laurent
polynomials $f \in \C[x^{\pm 1}, y^{\pm 1}, z^{\pm 1}]$, and write $\sfx^v
:= x^{a} y^{b} z^{c}$ where $v = (a,b,c) \in \Z^3$.

\begin{definition}\label{def:puzzle_pieces}
  Let $Q$ be a convex lattice polytope in $\Q^3$.
  \begin{enumerate}\itemsep=0pt
  \item We say that $Q$ is a \emph{length-one line segment} if
    there exist distinct points $v, w \in \Z^3$ such that the vertices of
    $Q$ are $v$, $w$ (in some order) and $Q \cap \Z^3 = \{v, w\}$.
    In this case, we set
    \[
    f_Q = \sfx^v + \sfx^w.
    \]
  \item We say that $Q$ is an \emph{$A_n$ triangle} if there
    exist distinct points $u, v, w \in \Z^3$ such that the vertices of $Q$
    are $u$, $v$, $w$ (in some order) and $Q \cap \Z^3 = \{u,
    v_0,v_1,\ldots,v_n\}$, where $v_0 = v$, $v_n = w$, and
    $v_0,v_1,\ldots,v_n$ are consecutive lattice points on the line
    segment from $v$ to $w$. In this case, we set
    \[
    f_Q = \sfx^u + \sum_{k=0}^n{n \choose k} \sfx^{v_k}.
    \]
  \end{enumerate}
  Observe that, in each case, the Newton polytope of $f_Q$ is $Q$.
\end{definition}

\begin{definition}
A \emph{lattice Minkowski decomposition} of a lattice polytope $Q\subset\Q^m$ is a decomposition $Q=Q_1+\cdots+Q_r$ of $Q$ as a Minkowski sum of lattice polytopes $Q_i$, $i\in\{1,2,\ldots,r\}$, such that the af\/f\/ine lattice generated by $Q\cap\Z^m$ is equal to the sum of the af\/f\/ine lattices generated by $Q_i\cap\Z^m$, $i\in\{1,2,\ldots,r\}$.
\end{definition}

\begin{example}
The Minkowski decomposition
\begin{center}
\begin{minipage}[c]{35pt}
	\includegraphics{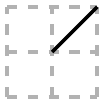}
\end{minipage}$+$
\begin{minipage}[c]{35pt}
	\includegraphics{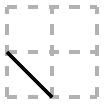}
\end{minipage}$=$
\begin{minipage}[c]{35pt}
	\includegraphics{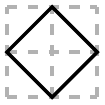}
\end{minipage}
\end{center}
is not a \emph{lattice} Minkowski decomposition, because the af\/f\/ine
lattices in the summands generate an index-two af\/f\/ine sublattice of
$\Z^2$.
\end{example}

\begin{definition}
  Let $Q$ be a lattice polytope in $\Q^m$. We say that lattice
  Minkowski decompositions
  \begin{gather*}
    Q = Q_1 + \cdots + Q_r, \qquad Q = Q'_1 + \cdots + Q'_s
  \end{gather*}
  of $Q$ are \emph{equivalent} if $r=s$ and if there exist $v_1,\ldots,v_r
  \in \Z^m$ and a permutation $\sigma$ of $\{1,2,\ldots,r\}$ such that
  $Q_i' = Q_{\sigma(i)} + v_i$.
\end{definition}

\begin{definition} Let $Q$ be a lattice polytope in $\Q^3$. A lattice
  Minkowski decomposition $Q = Q_1 + \cdots + Q_r$ of $Q$ is called
  \emph{admissible} if each $Q_i$ is either a length-one line segment
  or an $A_n$~triangle. Given an admissible Minkowski decomposition $Q
  = Q_1 + \cdots + Q_r$ of $Q$, we set
  \[
  f_{Q : Q_1,\ldots,Q_r} = \prod_{i=1}^{r} f_{Q_i}.
  \]
  Observe that the Newton polytope of $f_{Q : Q_1,\ldots,Q_r}$ is
  $Q$. Observe also that if
  \begin{gather*}
    Q = Q_1 + \cdots + Q_r, \qquad Q = Q'_1 + \cdots + Q'_s
  \end{gather*}
  are equivalent lattice Minkowski decompositions of $Q$, then $f_{Q
    : Q_1,\ldots,Q_r} = f_{Q : Q'_1,\ldots,Q'_s}$.
\end{definition}

\begin{definition}
  Let $P$ be a three-dimensional ref\/lexive polytope in $\Q^3$. Let
  $f$ be a Laurent polynomial in three variables. We say that $f$ is
  a \emph{three-dimensional Minkowski polynomial with Newton polytope
    $P$} if:
  \begin{enumerate}\itemsep=0pt
  \item $\Newt{f} = P$, so that
    \[
    f = \sum_{v \in P \cap \Z^3} a_v \sfx^v
    \]
    for some coef\/f\/icients $a_v$.
  \item $a_{\orig} = 0$.
  \item For each facet $Q\in\F{P}$ of $P$, we have
    \[
    \sum_{v \in Q \cap \Z^3} a_v \sfx^v = f_{Q:Q_1,\ldots,Q_r}
    \]
    for some admissible lattice Minkowski decomposition $Q = Q_1 +
    \cdots + Q_r$ of $Q$.
  \end{enumerate}
\end{definition}

Given a three-dimensional ref\/lexive polytope $P$, there may be:
\begin{enumerate}\itemsep=0pt
\item[1)] no Minkowski polynomial with Newton polytope $P$, if some facet
  of $P$ has no admissible lattice Minkowski decomposition;
\item[2)] a unique Minkowski polynomial with Newton polytope $P$, if there
  is a unique equivalence class of admissible lattice Minkowski
  decomposition for each facet of $P$;
\item[3)] many Minkowski polynomials with Newton polytope $P$, if there is
  some facet of $P$ with several inequivalent admissible lattice
  Minkowski decompositions.
\end{enumerate}
Up to isomorphism, there are $4319$ three-dimensional ref\/lexive
polytopes~\cite{KS98b}; $1294$ of them support no Minkowski polynomial,
and the remaining $3025$ together support $3747$
distinct\footnote{``Distinct'' here means ``distinct up to
  $\GL{3}{\Z}$-equivalence'', see
  Def\/inition~\ref{def:SL3_equivalence}.}  Minkowski polynomials.

\begin{figure}[!ht]
  \centering
  \includegraphics[scale=0.3]{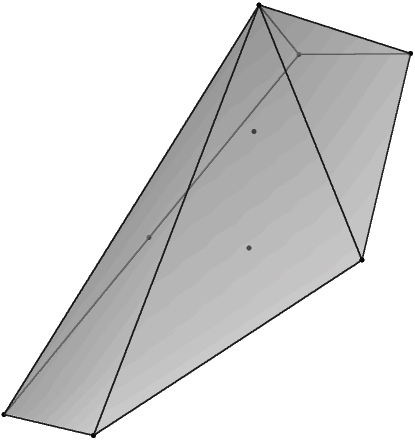}
  \caption{The polytope $P$ from Example~\ref{exa:122}.}
  \label{fig:122}
\end{figure}
\begin{example}\label{exa:122}
  Consider the three-dimensional ref\/lexive polytope $P$ with vertices:
  \[
  \{(-1, -1, -3), (1, 0, 0), (0, 1, 0), (0, 0, 1), (0, -1, -2), (-1,1, -1)\}.
  \]
  As facets, $P$ has four $A_1$ triangles, one $A_2$ triangle, and a
  pentagon. The pentagonal facet of $P$ has two inequivalent lattice
  Minkowski decompositions
\begin{center}
\begin{minipage}[c]{35pt}
	\includegraphics{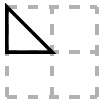}
\end{minipage}$+$
\begin{minipage}[c]{35pt}
	\includegraphics{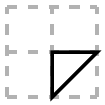}
\end{minipage}$=$
\begin{minipage}[c]{35pt}
	\includegraphics{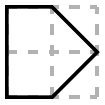}
\end{minipage}\hspace{2cm}
\begin{minipage}[c]{35pt}
	\includegraphics{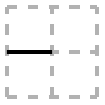}
\end{minipage}$+$
\begin{minipage}[c]{35pt}
	\includegraphics{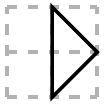}
\end{minipage}$=$
\begin{minipage}[c]{35pt}
	\includegraphics{F}
\end{minipage}
\end{center}
and thus there are two distinct Minkowski polynomials supported on $P$
  \begin{gather*}
    f_1  = x+y+z+\frac{1}{yz^2}+\frac{y}{xz}+\frac{2}{xz^2}+\frac{1}{xyz^3}+\frac{2}{z},\\
    f_2  = x+y+z+\frac{1}{yz^2}+\frac{y}{xz}+\frac{2}{xz^2}+\frac{1}{xyz^3}+\frac{3}{z}.
  \end{gather*}
\end{example}

\begin{remark}
  Recall from~Section~\ref{sec:introduction} that we denote by $N_f$ the lattice
  generated by the exponents of monomials of a Laurent polynomial~$f$.
  For a three-dimensional ref\/lexive polytope~$P$, the lattice
  generated by $P \cap \Z^3$ is equal to~$\Z^3$; thus whenever~$f$ is
  a three-dimensional Minkowski polynomial we have $N_f = \Z^3$. In
  dimension four there are precisely $16$ ref\/lexive polytopes~$P$ such
  that the lattice generated by $P \cap \Z^4$ is a proper sublattice
  of $\Z^4$~\cite{BK06}.
\end{remark}

\begin{remark}
  If $f$ is a Minkowski polynomial and $g$ is a Laurent polynomial
  related to $f$ via a~mutation $\varphi$ then in general $g$ will not
  be a Minkowski polynomial; indeed $\Newt{g}$ will in general not be
  ref\/lexive or even canonical
  (cf.~Corollary~\ref{cor:dim3_reflexive}).
\end{remark}

\section{Mutations between Minkowski polynomials}\label{sec:Minkowski_mutations}

Consider now the set of all three-dimensional Minkowski polynomials,
and partition this into smaller sets, which we call \emph{buckets},
such that $f$ and $g$ are in the same bucket if and only if the f\/irst
eight terms of the period sequences for $f$ and $g$ agree. There are
165 buckets. Using Proposition~\ref{prop:finitely_many} and a
computer search, we have determined all mutations of elements
of each bucket. For all except two buckets any two Laurent
polynomials $f$, $g$ in the bucket are connected by a sequence of
mutations that involves only Minkowski polynomials (and hence involves
only Laurent polynomials from that bucket). In fact:
\begin{enumerate}\itemsep=0pt
\item For $156$ of the buckets, any two Laurent polynomials $f$, $g$ in
  the bucket are connected by a~sequence of \emph{width two} mutations
  that involves only Minkowski polynomials in that bucket.
\item For $7$ of the buckets, any two Laurent polynomials $f$, $g$ in
  the bucket are connected by a~sequence of \emph{width two and width three}
  mutations that involves only Minkowski polynomials in that bucket.
\end{enumerate}
For the remaining two buckets\footnote{Numbering as in Appendix~\href{http://www.emis.de/journals/SIGMA/2012/094/sigma12-094Appendices.pdf}{A},
  these are the buckets with IDs~148 and~161.}, any two Laurent
polynomials $f$, $g$ in the bucket are connected by a sequence of
mutations, but it is not possible to insist that all of the Laurent
polynomials involved are Minkowski polynomials. (It is, however,
possible to insist that all of the Laurent polynomials have ref\/lexive
Newton polytope.)

The sequence of mutations connecting two Minkowski polynomials $f$ and
$g$ is far from unique, but representative such mutations are shown in
Appendix~\href{http://www.emis.de/journals/SIGMA/2012/094/sigma12-094Appendices.pdf}{B}. The mutations shown there
suf\/f\/ice to connect any two three-dimensional Minkowski polynomials in
the same bucket.

\begin{corollary}
  Three-dimensional Minkowski polynomials $f$ and $g$ have the same
  period sequence if and only if they have the same first eight terms
  of the period sequence.
\end{corollary}

\begin{proof}
  Any two Minkowski polynomials in the same bucket are connected by a
  sequence of mutations.
\end{proof}

\begin{proposition}
  If $f$ and $g$ are three-dimensional Minkowski polynomials with the
  same period sequence, then there exist a flat family $\pi :
  \mathcal{X} \to \Sigma$ over a possibly-reducible rational curve
  $\Sigma$ and two distinct points $0$,~$\infty \in \Sigma$ such that
  $\pi^{-1}(0)$ is isomorphic to the toric variety $X_f$
  and $\pi^{-1}(\infty)$ is isomorphic to the toric variety $X_g$.
\end{proposition}

\begin{proof}
  Combine the results in Appendix~\href{http://www.emis.de/journals/SIGMA/2012/094/sigma12-094Appendices.pdf}{B} with~\cite[Theorem~1.3]{I12}.
\end{proof}

In particular this gives a geometric proof that if $f$ and $g$ are
three-dimensional Minkowski polynomials with the same period sequence,
then $X_f$ and $X_g$ have the same Hilbert series.  One might hope
that substantially more is true: for example that given a bucket of
Minkowski polynomials $\{f_i : i \in I\}$ corresponding to a Fano
manifold $X$, the toric varieties $X_{f_i}$, $i \in I$, all lie in the
same component of the Hilbert scheme, and that $X$ lies in this
component too.  It is possible that this hope is too na\"\i ve, as the
geometry of Hilbert schemes is in many ways quite pathological, but in
any case, as things currently stand, results of this form seem to be
out of reach.

\begin{remark}
  One might regard our results as further evidence that the mirror to
  a Fano manifold~$X$ should be some sort of cluster variety $X^\vee$
  together with a function (the superpotential) on~$X^\vee$.  This is
  familiar from the work of Auroux~\cite{Auroux}.  Given a special
  Lagrangian cycle and a complex structure on a Fano manifold~$X$,
  mirror symmetry associates to this choice a Laurent polynomial.  The
  Laurent polynomial depends on our choice: there is a
  wall-and-chamber structure on the space of parameters, and Laurent
  polynomials from neighbouring chambers are expected to be related by
  an elementary birational transformation.  Reformulating this: the
  mirror to~$X$ should be a pair~$(X^\vee,W)$ where~$X^\vee$ is
  obtained by gluing tori (one for each chamber) along birational
  transformations (coming from wall-crossing), and the superpotential~$W$ is a global function on~$X^\vee$.  When $X$ is two-dimensional,
  the elementary birational transformations that occur here are
  closely related to those in the theory of cluster algebras~\cite{FZ}.  In our three-dimensional situation, we need to allow a
  still more general notion of cluster algebra.

\looseness=1
  More precisely: the fundamental objects in the theory of cluster
  algebras are \emph{seeds} and their mutations.  A seed is a
  collection of combinatorial and algebraic data: the so-called
  exchange matrix and a basis for a f\/ield $k = \Q(x_1,\dots,x_n)$ of
  rational functions of $n$ variables, called a cluster.  Seeds can be
  transformed into new seeds through a process called
  \emph{mutation}. The \emph{Laurent phenomenon} theorem of~\cite{FZ}
  says that any cluster variable is a Laurent polynomial when
  expressed in terms of any other cluster variable in any other
  mutation equivalent seed. This produces a large subalgebra in $k$,
  called the \emph{upper bound}, consisting of those elements that are
  Laurent polynomials when expressed in terms of cluster variables of
  any seed.  The ``Laurent phenomenon'' arising from mirror symmetry
  for Fano manifolds (i.e.~the presence of a Laurent polynomial that
  remains a Laurent polynomial under some collection of ele\-men\-tary
  birational transformations) does not f\/it into the framework of
  cluster algebras. In cluster algebras the number of ways to mutate a
  seed naturally is always less than or equal to the transcendence
  degree of $k$, whereas the Laurent polynomials that are mirror dual
  to a del~Pezzo surface $S$ have $\chi(S) = \rho(S)+2$ ways to mutate
  them~\cite{CMG12, GU10} and the transcendence degree of $k$ here is
  $2$.  Furthermore, in contrast to the case of cluster algebras, the
  transcendence degree of the upper bound in this setting equals just
  one; in fact, the upper bound is just the ring of polynomials
  generated by $W$.  In dimension three there are additional
  complications: in the Fano setting the exchange polynomials need not
  be binomials (cf.~\cite{LP}), and we do not see how to def\/ine the
  2-form that occurs in the cluster theory.  In this paper we do not
  consider the problem of def\/ining seeds in dimension three, nor do we
  single out the relevant notion of mutation of seeds. These are two
  important questions for further study.
\end{remark}

\subsection*{Acknowledgements}

We thank Alessio Corti and Vasily Golyshev for many useful conversations, the referees for perceptive and helpful comments, John Cannon for providing copies of the computer algebra software {\magma}, and Andy Thomas for technical assistance. This research is supported by a~Royal Society University Research Fellowship; ERC Starting Investigator Grant number~240123; the Leverhulme Trust; Kavli Institute for the Physics and Mathematics of the Universe (WPI); World Premier International Research Center Initiative (WPI Initiative), MEXT, Japan; Grant-in-Aid for Scientif\/ic Research (10554503) from Japan Society for Promotion of Science and Grant of Leading Scientif\/ic Schools (N.Sh. 4713.2010.1); and EPSRC grant EP/I008128/1.

\pdfbookmark[1]{References}{ref}
\LastPageEnding
\end{document}